\documentclass[12pt]{article}

\usepackage[latin1]{inputenc}
\usepackage[english]{babel}
\usepackage{amssymb,amsmath,amscd,amstext, amsthm}
\usepackage{amsfonts}
\usepackage[dvips]{graphicx}
\usepackage{psfrag}
\usepackage{epsfig}
\usepackage{palatino}
\usepackage{xypic}

\usepackage{graphicx}
\usepackage[all]{xy}
\usepackage{amscd}

\newtheorem{thm}{Theorem}[section]    % Standard theorem environment
\newtheorem{lem}[thm]{Lemma}          % Lemma environment with numbering
\theoremstyle{definition}
\newcommand {\Z} {\mathbb {Z}}
\newcommand {\Q} {\mathbb {Q}}
\newcommand {\ZG} {\mathbb {Z}[G]}

\newcommand {\Rk} {{\rm rk}_{\Z}}

\renewcommand{\vec}[1]{\boldsymbol{#1}}
\newcounter{eno}
\newcounter{ch}
\newcounter{sec}
\newcounter{no}
\def \sec#1 {\stepcounter{sec} \bigskip\par{\underline{\bf \S \arabic{sec} {\bf #1}  \setcounter{eno}{1} \setcounter{no}{1}}}}
\def \prop#1 {\par{\bf Proposition \arabic{sec}.\arabic{no}} \stepcounter{no} {\it #1}}
\def \cor#1 {\par{\bf Corollary \arabic{sec}.\arabic{no}} \stepcounter{no} {\it #1}}
\def \lem#1 {\par{\bf Lemma \arabic{sec}.\arabic{no}} \stepcounter{no} {\it #1}}
\def \thm#1 {\par{\bf Theorem \arabic{sec}.\arabic{no}} \stepcounter{no} {\it #1}}
\def \conj#1 {\par{\bf Conjecture \arabic{sec}.\arabic{no}} \stepcounter{no} {\it #1}}
\def \define#1 {\par{\bf Definition \arabic{sec}.\arabic{no} } \stepcounter{no} {\it #1} \hspace{1mm}}

\def\co{\colon\thinspace}
\def \proof#1 {\bigskip \hspace{2mm}Proof:$\,\,\,$  #1 ${}$\hfill $\Box$
\,\, \newline }

\begin{document}

\begin{center}
{\bf \huge Quillen's plus construction and the D(2) problem}

{W. H. Mannan}

\end{center}

\bigskip
{\bf MSC}: 57M20, 19D06, 57M05 \hfill {\bf Keywords}: {D2 problem}, {Quillen plus}
\hfill construction

\begin{quote}
Given a finite connected 3--complex with cohomological dimension 2, we show it may be constructed up to homotopy  by applying the Quillen plus construction to the Cayley complex of a finite group presentation.  This reduces the D(2) problem to a question about perfect normal subgroups.
\end{quote}

\sec{Introduction}

Given a finite cell complex one may ask what the minimal dimension of a finite cell complex in its homotopy type is. If $n\neq 2$ and the cell complex has cohomological dimension $n$ (with respect to all coefficient bundles), then the cell complex is in fact homotopy equivalent to a finite $n$--complex (a cell complex whose cells have
 dimension at most $n$).  Although this has been known for around forty years (for $n>2$ it is proved in \cite{Wall} and for $n=1$ it follows from \cite{Swan}, \cite{Stal}), it is an open question whether or not this holds when $n=2$.  This question is known as Wall's D(2) problem:

{ \bf Let $X$ be a finite 3--complex with $H^3(X;\beta)=0$ for all coefficient bundles $\beta$.  Must $X$ be homotopy equivalent to a finite 2--complex?
}

If $X$ (as above) is not homotopy equivalent to a finite 2--complex, we say it is a counterexample which solves the D(2) problem.

For connected $X$ with certain fundamental groups, it has shown been shown that $X$ must be homotopy equivalent to a finite 2--complex (see for example  \cite{John}, \cite{Edwa}, \cite{Mann}).  However no general method has been forthcoming.

Also, whilst potential candidates for counterexamples have been constructed (see \cite{Beyl}, \cite{Brid} ), no successful method has yet emerged for verifying that they are not homotopy equivalent to finite 2--complexes.

In \S2 we apply the Quillen plus construction to  connected 2--complexes, resulting in cohomologically 2--dimensional 3--complexes.  These are therefore candidates for counterexamples which solve the D(2) problem.  In \S3 we show that in fact all finite connected cohomologically 2--dimensional 3--complexes arise this way, up to homotopy equivalence.

Finally, in \S4 we use these results to reduce the D(2) problem to a question about perfect normal subgroups.  This allows us to generalize existing approaches to the D(2) problem such as \cite[Theorem I]{John1} and \cite[Theorem 3.5]{Harl}.

\bigskip
Before moving on to the main argument we make a few notational points.  All modules are right modules except where a left action is explicitly stated.  The basepoint of a Cayley complex is always assumed to be its 0--cell.

If $X$ is a connected cell complex with basepoint, we denote its universal cover $\tilde{X}$.  Given two based loops $\gamma_1, \gamma_2 \in \pi_1(X)$ their product $\gamma_1\gamma_2$ is the composition whose initial segment is $\gamma_2$ and final segment is $\gamma_1$.   With this convention, we have a natural right action of $\pi_1(X)$ on the cells of $\tilde{X}$.  Let $G=\pi_1(X)$. We can regard the associated chain complex of $\tilde{X}$ as an algebraic complex of right modules over $\ZG$.  We follow \cite{John1} in denoting this algebraic complex $C_*(X)$.  Note that this differs from the convention in other texts.  Thus in particular $C_*(X)$ and $C_*(\tilde{X})$ have the same underlying sequence of abelian groups, but the former is a sequence of modules over $\ZG$ whilst the latter is a sequence of modules over $\Z[\pi_1(\tilde{X})]=\Z$.

If $Y$ is a subcomplex of $X$ then $C_*(Y)$ is a sequence of right
modules over $\pi_1(Y)$.  Let $E=\pi_1(Y)$.  The induced map $E
\to G$ yields a left action of $E$ on $\ZG$.  Thus we have an
algebraic complex $C_*(Y) \otimes_E \ZG$ over $\ZG$.  The
inclusion $Y \subset X$ induces a chain map $C_*(Y) \otimes_E \ZG
\longrightarrow C_*(X)$.  The complex $C_*(X,Y)$ is defined to be
the relative chain complex associated to this chain map.

The basepoint allows us to interchange between coefficient bundles over $X$ and right modules over $\ZG$.  Thus for a right module $N$ we have:
$$H^n(X;N)=H^n(C_*(X); N)$$

A left module over $\ZG$ may be regarded as a right module over $\ZG$, where right multiplication by a group element is defined to be left multiplication by its inverse.  Hence a left module $M$ may also be regarded as a coefficient bundle and we have: $$H_n(X;M)=H_n(C_*(X); M),\qquad \qquad H_n(X,Y;M)=H_n(C_*(X,Y) ; M)$$

Given a finitely generated Abelian group $A$ we may regard it as a finitely generated module over $\Z$.  Thus $A \otimes_\Z \Q$ is a finite dimensional vector space over $\Q$.   The dimension of this vector space will be denoted $\Rk(A)$.

Finally given a group $G$ and elements $g,h \in G$, we follow the convention that $[g,h]$ denotes the element $ghg^{-1}h^{-1}$.

\sec{The plus construction applied to a Cayley complex}

Let $\varepsilon= \langle g_1, \cdots,g_n \,|\, R_1, \cdots, R_m \rangle$ be a finite presentation for a group $E$.  We say a normal subgroup of $E$ is {\sl finitely closed} when it is the normal closure in $E$ of a finitely generated subgroup.  Let $K \lhd E$ be finitely closed and perfect (so $K=[K, K]$).  Let $\mathcal{K}_\varepsilon$ denote the Cayley complex associated to $\varepsilon$.

\thm{{\rm (Quillen, see \cite[Theorem 5.2.2]{Rose})} There is a 3--complex $\mathcal{K}_\varepsilon^+$,  containing $\mathcal{K}_\varepsilon$ as a subcomplex, such that the inclusion $\mathcal{K}_\varepsilon \hookrightarrow \mathcal{K}_\varepsilon^+$ induces the quotient map $E \to E/K$ on fundamental groups and $H_*(\mathcal{K}_\varepsilon^+,\mathcal{K}_\varepsilon; M) =0$ for all left modules $M$ over $\Z[E/K]$.  Further, given another such 3--complex $X$, there is a homotopy equivalence $\mathcal{K}_\varepsilon^+ \to X$ extending the
identity map of the common subspace $\mathcal{K}_\varepsilon$.}

In fact we may construct $\mathcal{K}_\varepsilon^+$ explicitly, using the fact that $K$ is finitely closed to ensure that we end up with a finite cell complex.  Let $k_1,\cdots k_r \in K$ generate a subgroup of $E$ whose normal closure (in $E$) is $K$.  As $K=[K,K]$, each $k_i$ may be expressed as a product of commutators: $k_i = \prod_{j=1}^{m_i} [a_{ij}, b_{ij}]$ with each $a_{ij}, b_{ij} \in K$.  Then each $a_{ij}, b_{ij}$ may be represented by words $A_{ij}, B_{ij}$ in the $g_l$, $l=1, \cdots,n$.  For each $i=1,\cdots,r$ attach a 2--cell $E_i$ to $\mathcal{K}_\varepsilon$ whose boundary corresponds to the word $\prod_{j=1}^{m_i} [A_{ij}, B_{ij}]$.  Denote the resulting chain complex $\mathcal{K}_\varepsilon'$.

The chain complex $C_*(\mathcal{K}_\varepsilon)$ may written: $$C_*(\mathcal{K}_\varepsilon) \colon C_2(\mathcal{K}_\varepsilon) \stackrel{\partial_2}\to C_1(\mathcal{K}_\varepsilon) \stackrel{\partial_1}\to C_0(\mathcal{K}_\varepsilon)
$$

Note that the boundary map $\partial_2$ applied to a 2--cell is the Fox free differential: $\partial \co F_{\{g_1,\cdots,g_n\}} \to C_1(\mathcal{K}_\varepsilon) $, applied to the word which the 2--cell bounds (see \cite[\S48]{John1} and \cite{Fox}).  Let $\vec{e_i}$ denote the generator in $C_1(\mathcal{K}_\varepsilon)$
representing the generator $g_i$.  The free Fox differential is then characterized by:

i) $\partial g_i=\vec{e_i}$ for all $i=1,\cdots, n$

ii) $\partial(AB)=\partial(A)B+\partial(B)$ for all words $A,B$.

\bigskip
Clearly the inclusion $\mathcal{K}_\varepsilon \hookrightarrow \mathcal{K}_\varepsilon'$ induces the quotient map $E \to E/K$ on fundamental groups.   There is a right action of $\Z[E/K]$ on itself.  Further there is a left action of $E$ on $\Z[E/K]$.

\lem{As an algebraic complex of right $\Z[E/K]$ modules $C_*(\mathcal{K}_\varepsilon')$ may be written:  $$C_*(\mathcal{K}_\varepsilon')\colon  C_2(\mathcal{K}_\varepsilon)\otimes_E \Z[E/K] \oplus \Z[E/K]^r \stackrel{\partial_2 \oplus 0}\to C_1(\mathcal{K}_\varepsilon)\otimes_E \Z[E/K] \stackrel{\partial_1 }\to C_0(\mathcal{K}_\varepsilon)\otimes_E \Z[E/K]
$$}

\hspace{2mm}Proof:  The boundary of $E_i$ is given by the free Fox differential $\partial$, applied to the word $\prod_{j=1}^{m_i} [A_{ij}, B_{ij}]$.  However, $$\partial \prod_{j=1}^{m_i} [A_{ij}, B_{ij}] = \sum_{j=1}^{m_i} [\partial A_{ij} + \partial B_{ij} - \partial A_{ij} -\partial B_{ij}]=0$$
as each $A_{ij},B_{ij}$ represents an element of $K$ and hence is trivial in $\pi_1(\mathcal{K}_\varepsilon')=E/K$. \hfill $\Box$

\bigskip
Each $E_i$ therefore generates an element of $H_2(\tilde{\mathcal{K}_\varepsilon'};\Z)$.  By the Hurewicz isomorphism Theorem we have isomorphisms $H_2(\tilde{\mathcal{K}_\varepsilon'};\Z) \cong \pi_2(\tilde{\mathcal{K}_\varepsilon'}) \cong \pi_2(\mathcal{K}_\varepsilon')$ coming from the Hurewicz homomorphism and the covering map respectively.  Let $\psi_i \co S^2 \to \mathcal{K}_\varepsilon'$ represent the element of $\pi_2(\mathcal{K}_\varepsilon')$ which corresponds to $E_i$ under these isomorphisms.

For each $i \in 1, \cdots r$ we then attach a 3--cell $B_i$ to $\mathcal{K}_\varepsilon'$ via the attaching map:\newline
$\psi_i\co\partial B_i \to \mathcal{K}_\varepsilon'$.  Let $\mathcal{K}_\varepsilon''$ denote the resulting 3--complex.  Then we have  $C_*(\mathcal{K}_\varepsilon''):$
$$
\Z[E/K]^r \stackrel{\partial_3}\to C_2(\mathcal{K}_\varepsilon)\otimes_E \Z[E/K] \oplus \Z[E/K]^r \stackrel{\partial_2 \oplus 0}\to C_2(\mathcal{K}_\varepsilon)\otimes_E \Z[E/K] \stackrel{\partial_1 }\to C_0(\mathcal{K}_\varepsilon)\otimes_E \Z[E/K]
$$
where $\partial_3$ is inclusion of the second summand.

Hence we have:

\lem{$H_*(\mathcal{K}_\varepsilon'',\mathcal{K}_\varepsilon; M) =0$ for all left modules $M$ over $\Z[E/K]$.}

\proof{We have  the following relative complex: $$C_*(\mathcal{K}_\varepsilon'',\mathcal{K}_\varepsilon) \colon \quad \Z[E/K]^r \tilde{\longrightarrow} \Z[E/K]^r \to 0 \to 0$$}

Thus by Theorem 2.1 we may conclude that $\mathcal{K}_\varepsilon''$ has the homotopy type of $\mathcal{K}_\varepsilon^+$.

\lem{The complex $\mathcal{K}_\varepsilon''$ is cohomologically
2--dimensional.}

\proof{The inclusion $\iota\co \mathcal{K}_\varepsilon
\hookrightarrow \mathcal{K}_\varepsilon''$ induces a chain
homotopy equivalence:
$$C_*(\mathcal{K}_\varepsilon)\otimes_E \Z[E/K] \to
C_*(\mathcal{K}_\varepsilon'')$$}

\cor{We may choose $\mathcal{K}_\varepsilon^+$ to be the
cohomologically 2--dimensional finite 3--complex
$\mathcal{K}_\varepsilon''$.}

\sec{Cohomologically 2--dimensional 3--complexes}

 Let $X$ be a finite connected 3--complex with $H^3(X;\beta)=0$ for all coefficient bundles $\beta$.  In this section we will show that up to homotopy, $X$ arises as the Quillen plus construction applied to a finite Cayley complex.

Let $T$ be a maximal tree in the 1--skeleton of $X$.  The quotient map $X \to X/T$ is a homotopy equivalence.   Hence we may assume without loss of generality that $X$ has one $0$--cell.  We take this to be the basepoint of $X$ and any complexes obtained from $X$ by adding or removing cells.  Also we set $G=\pi_1(X)$ with respect to this basepoint.

 Let $C_*(X)$ be denoted by:

 $$
 F_3 \stackrel{\partial_3}\to F_2 \stackrel{\partial_2}\to F_1 \stackrel{\partial_1}\to F_0
$$
where the $F_i$, $i=0,1,2,3$, are free modules over $\ZG$ and the $\partial_i$ are linear maps over $\ZG$.

We have $H^3(X;F_3) =0$ so in particular there exists $\phi$ such the following diagram commutes:

\xymatrix{
&&&& F_3 \ar[d]^1\ar[r]^{\partial_3}& F_2\ar[dl]^\phi \ar[r]^{\partial_2}& F_1 \ar[r]^{\partial_1}&F_0\\
&&&& F_3
}

\bigskip
Hence $\partial_3$ is the inclusion of the first summand $\partial_3 \co F_3 \hookrightarrow \partial_3(F_3) \oplus S = F_2$, where $S$ is the kernel of $\phi$.  Let $X'$ denote the wedge of $X$ with one disk for each 3--cell in $X$.  Then the inclusion of cell complexes $X \hookrightarrow X'$ is a homotopy equivalence and:

$$
C_*(X') \co \qquad F_3 \stackrel{\partial_3'}\to F_2 \oplus F_3' \stackrel{\partial_2'}\to F_1 \oplus F_3' \stackrel{\partial_1'}\to F_0
$$

Here $F_3' \cong F_3$ and the maps are defined as follows:

$\partial_1'$ restricts to $\partial_1$ on $F_1$ and restricts to $0$ on $F_3'$,

$\partial_2'= \left(\begin{array}{cc}\partial_2&0\\0&1\end{array}\right)$

and $\partial_3'$ is $\partial_3 \co F_3 \to F_2$ composed with the natural inclusion: $F_2 \hookrightarrow F_2 \oplus F_3'$.  Thus $\partial_3'$ is the inclusion into the first summand: $\partial_3'\co F_3 \hookrightarrow \partial_3'F_3 \oplus S \oplus F_3'$.

Let $m$ denote the number of 2--cells in $X$.  The submodule $S \oplus F_3' \subset (\partial_3'F_3 \oplus S) \oplus F_3' $ is isomorphic to $S \oplus F_3 \cong F_2$ and hence has a basis $\vec{x_1},\cdots, \vec{x_m} \in F_2 \oplus F_3'$.

The cell complex $X'$ has one 0--cell, so $F_0 \cong \ZG$.  Let $n$ denote the number of 1--cells in $X'$.  Then each 1--cell corresponds to a generator $g_i,\, i \in [1,\cdots,n]$ of $G$.  Let $\{\vec{e_1},\cdots,\vec{e_n}\}$ form the corresponding basis for $F_1 \oplus F_3'$.

Let $r$ denote the number of 2--cells in $X'$.  The attaching map for each 2--cell maps the boundary of a disk round a word in the $g_i$.  For each 2--cell let $R_j,\, j\in [1,\cdots,r]$ denote this word.  Let $\{\vec{E_1}, \cdots,\vec{E_r}\}$ form the corresponding basis for $F_2 \oplus F_3'$.  Thus we have a presentation $G = \langle g_1,\cdots,g_n \,|\,R_1,\cdots,R_r\rangle$.

We may therefore express each $\vec{x_i}$ as a linear combination of the $\vec{E_j}$.  Thus for some integers $v_i$ and sequences $j_{i1},\cdots,j_{iv_i} \in \{1,\cdots,r\}$ we have:  $$\vec{x_i}= \sum_{l=1}^{v_i} \vec{E_{j_{il}}}\lambda_{il}\sigma_{il}$$ with each $\lambda_{il} \in G$ and $\sigma_{il} \in \{1,-1\}$.  For each $i \in [1,\cdots, m]$, $l\in[1,\cdots, v_i]$ let $w_{il}$ be a word in the $g_k,\, k=1,\cdots,n$, representing $\lambda_{il}$.  Now for each $i=1,\cdots m$, let: $$S_i=
\prod_{l=1}^{v_i} w_{il}^{-1}R_{j_{il}}^{\sigma_{il}} w_{il}$$

For each $i \in \{1,\cdots m\}$, attach a 2--cell $a_i$ to $X'$ by mapping the boundary of the disk around the path in the 1--skeleton of $X'$ corresponding to the word $S_i$.  Let $Z$ denote the resulting finite cell complex.  Note that each word $S_i$ corresponds to a trivial element of $G$, so the inclusion $X' \subset Z$ induces an isomorphism $\pi_1(X') \cong \pi_1(Z)$.  Hence we may write $C_*(Z) \colon$

$$ C_*(Z) \colon \quad
F_3 \stackrel{\partial_3''}\to (F_2 \oplus F_3')\oplus F_2' \stackrel{\left(\begin{array}{cc}\partial_2'&\partial_2''\end{array}\right) } \longrightarrow (F_1 \oplus F_3') \stackrel{\partial_1'}\to F_0
$$

where $\partial_3''$ is understood to be $\partial_3'\co F_3 \to (F_2 \oplus F_3')$ composed with the natural inclusion $(F_2 \oplus F_3') \hookrightarrow (F_2 \oplus F_3')\oplus F_2'$.

For $i=1,\cdots,m$ let $\vec{A_i}$ be the basis element of $F_2'$ corresponding to the 2--cell $a_i$.  Recall the Fox free differential, $\partial$. We have:$$
\partial_2'' \vec{A_i} = \partial S_i = \sum_{l=1}^{v_i} \partial (w_{il}^{-1}R_{j_{il}}^{\sigma_{il}} w_{il}) = \sum_{l=1}^{v_i} \partial_2' \vec{E_{j_{il}}} \lambda_{il} \sigma_{il} = \partial_2'\vec{x_i}
$$

Thus $\vec{A_i}-\vec{x_i}$ represents a class in $H_2(\widetilde{Z^{(2)}}; \Z)$ which is isomorphic to $\pi_2(Z^{(2)})$ via the Hurewicz isomorphism composed with the map $\pi_2(\widetilde{Z^{(2)}}) \to \pi_2(Z^{(2)})$ induced by the covering map.  Let $\psi_i \co S^2 \to {Z^{(2)}}$ represent the corresponding element of $\pi_2({Z^{(2)}})$.

Then for each $i=1,\cdots,m$ we may attach a 3--cell $b_i$ to $Z$ via the map $\psi_i$.  We denote the resulting complex $X''$.

\lem{The inclusion $\iota \co X' \subset X''$ is a homotopy equivalence.}

\vspace{-4mm}

\proof{Starting with $X'$, for each $i$ we attached a 2--cell $a_i$ with contractable boundary in $X'$, and then attached a 3--cell $b_i$ with $a_i$ as a free face.  Thus $X''$ is obtained from $X'$ through a series of cell expansions and the inclusion $X' \subset X''$ is a simple homotopy equivalence.}

\vspace{-4mm}

Let $Y$ denote the subcomplex of $X''$ consisting of the 1--skeleton, $X''^{(1)}$, together with the $a_i$, $i=1,\cdots,m$.  Let $\varepsilon$ denote the group presentation
$\langle g_1,\cdots,g_n \,|\,S_1,\cdots,S_m\rangle$ and let $E$ denote the underlying group.  By construction we have: $Y=\mathcal{K}_\varepsilon$.

Let $k_1, \cdots, k_r \in E$ denote the elements represented by the words $R_1,\cdots,R_r$.  Let $K$ denote the normal closure in $E$ of $k_1, \cdots,k_r$.  By construction then, $K$ is finitely closed and we have a short exact sequence of groups:$$1 \to K \to E \to G \to 1$$

\vspace{-2mm}

\lem{$K$ is a perfect group.}

\hspace{2mm} Proof:  {Clearly $\ZG$ is a right module over itself and there is a left action of $E$ on $\ZG$.  The algebraic complex $C_*(\mathcal{K}_\varepsilon) \otimes_E \ZG$ is given by:
$$
F_2' \stackrel{\partial_2''} \rightarrow (F_1 \oplus F_3') \stackrel{\partial_1'}\to F_0
$$
Now consider $C_*(X')$:
$$
 F_3 \stackrel{\partial_3'}\to F_2 \oplus F_3' \stackrel{\partial_2'}\to F_1 \oplus F_3' \stackrel{\partial_1'}\to F_0
$$
As $\tilde{X}'$ is simply connected, we have ${\rm ker }(\partial_1') = {\rm Im}(\partial_2')$.
\newline
\newline
Recall that $F_2 \oplus F_3' = \partial_3' (F_3) \oplus S \oplus F_3'$ and that $S \oplus F_3'$ has basis $\vec{x_1}, \cdots,\vec{x_m}$.  Clearly $\partial_2'$ restricts to $0$ on $\partial_3' (F_3)$, so ${\rm ker }(\partial_1') = {\rm Im}(\partial_2')$ which is generated by the $\partial_2'(\vec{x_i})$.
\newline
\newline
Also recall that $\partial_2'\vec{x_i}=\partial_2'' \vec{A_i}$.  Hence ${\rm ker }(\partial_1') = {\rm Im}(\partial_2'')$ and $H_1(C_*(\mathcal{K}_\varepsilon) \otimes_E \ZG)=0$.
\newline
\newline
However by restricting coefficients $C_*(\mathcal{K}_\varepsilon)$ may be regarded as an algebraic complex of free modules over $\Z[K]$. Hence we have:\newline \newline
$K/[K,K] = H_1(K;\Z) = H_1(C_*(\mathcal{K}_\varepsilon) \otimes_{K} \Z) =H_1(C_*(\mathcal{K}_\varepsilon) \otimes_E \ZG)=0$.
\newline\newline where $\Z$ is regarded as having a trivial left $K$--action.} \hfill $\Box$

\lem{$X''= \mathcal{K}_\varepsilon^+$ where $+$ is taken with respect to $K$.}

\proof{We may identify $\mathcal{K}_\varepsilon$ with the subcomplex $Y \subset X''$.  The inclusion $\ell \co \mathcal{K}_\varepsilon \hookrightarrow X''$ then induces the quotient map $E \to E/K$ on fundamental groups.  By Theorem 2.1 it is sufficient to show that $H_*(X'',Y;M)=0$ for all left  coefficient modules $M$.
\newline
Let $\ell_*\co C_*(\mathcal{K}_\varepsilon) \otimes_E \ZG \to C_*(X'')$ be the chain map induced by the inclusion $\ell\co \mathcal{K}_\varepsilon \hookrightarrow X''$. We have the following commutative diagram:
\newline
\begin{eqnarray*}
F_2' \,\,\,\quad \stackrel{\partial_2''}\longrightarrow \quad \,\,\,\,(F_1 \oplus F_3') \stackrel{\partial_1'}\to F_0\,\,\,\,\\
\downarrow^{\ell_2}\qquad\,\,\, \,\qquad\qquad\downarrow^{\ell_1}\,\,\,\,\qquad
\downarrow^{\ell_0}\\
F_3 \oplus F_2''\stackrel{\left(\begin{array}{cc}\partial_3''&\partial_3'''\end{array}\right)}\to (F_2 \oplus F_3')\oplus F_2' \stackrel{\left(\begin{array}{cc}\partial_2'&\partial_2''\end{array}\right) } \longrightarrow (F_1 \oplus F_3') \stackrel{\partial_1'}\to F_0\,\,\,\,
\end{eqnarray*}\newline
where $F_2''$ has a basis $\vec{D_1},\cdots,\vec{D_m}$ corresponding to the 3--cells $b_1,\cdots,b_m$, so for $i=1,\cdots,m$ we have $\partial_3'''(D_i) =\vec{A_i} -\vec{x_i}$.    Here $\ell_0$ and $\ell_1$ are the identity maps and $\ell_2$ is the inclusion of the second summand.
\newline
\newline
We have $(F_2 \oplus F_3') = \partial_3''F_3 \oplus (S \oplus F_3')$. Hence we have $(F_2 \oplus F_3')\oplus F_2' = \partial_3''F_3 \oplus (S \oplus F_3') \oplus F_2'$.
\newline\newline
 The submodule  $(S \oplus F_3')$ has basis $\vec{x_1},\cdots, \vec{x_m}$.  The submodule $F_2'$ has basis $\vec{A_1},\cdots,\vec{A_m}$.  Also $\partial_3'''F_2''$ has basis $\vec{A_1}-\vec{x_1},\cdots,\vec{A_m}-\vec{x_m}$.  Hence we have the following equality of submodules: $(S \oplus F_3') \oplus F_2' = \partial_3'''F_2'' \oplus F_2'$ \newline \newline
Thus: $$(F_2 \oplus F_3')\oplus F_2' = \partial_3''F_3 \oplus \partial_3'''F_2'' \oplus F_2'$$   The relative chain complex $C_*(X'',Y)$ is therefore given by:
$$
F_3 \oplus F_2'' {\tilde{\longrightarrow}} \partial_3''F_3 \oplus \partial_3'''F_2'' \longrightarrow 0 \longrightarrow 0
$$ and $H_*(X'',Y;M)=0$ for all left coefficient modules $M$ as required.
}

As $X \sim X''$, we have proved the following theorem:

\thm{Let $X$ be a finite connected 3--complex with $H^3(X;\beta)=0$ for all coefficient bundles $\beta$.  Then $X$ has the homotopy type of $\mathcal{K}_\varepsilon^+$ for some finite presentation $\varepsilon$ of a group $E$, where $+$ is taken with respect to some perfect finitely closed normal subgroup $K \lhd E$. }

\sec{Implications for the D(2) problem}

The D(2) problem asks if every finite cohomologically 2--dimensional 3--complex must be homotopy equivalent to a finite 2--complex.  Clearly a counterexample must have a connected component which is also cohomologically 2--dimensional but not homotopy equivalent to a finite 2--complex.  By Theorem 3.4 this component must have the homotopy type of $\mathcal{K}_\varepsilon^+$ for some finite presentation $\varepsilon$ of a group $E$, where $+$ is taken with respect to some perfect finitely closed normal subgroup $K \vartriangleleft E$.

Conversely, by Corollary 2.5, given any finite presentation $\varepsilon$ of a group $E$ together with some perfect finitely closed normal subgroup $K \vartriangleleft E$ we have a cohomologically 2--dimensional finite 3--complex, $\mathcal{K}_\varepsilon^+$.  It follows that the D(2) problem is equivalent to:

{\bf Given a finite presentation $\epsilon$ for a group $E$, and a finitely closed perfect normal subgroup $K \vartriangleleft E$, must  $\mathcal{K}_\varepsilon^+$ be homotopy equivalent to a finite 2--complex?}

Suppose that we have a homotopy equivalence $\mathcal{K}_\varepsilon^+ \sim Y$ for some finite 2--complex $Y$.  Let $T$ be a maximal tree in the $1$--skeleton of $Y$.  The quotient map $Y \to Y/T$ is a homotopy equivalence so $Y \sim \mathcal{K}_\mathcal{G}$ for some finite presentation $\mathcal{G}$ of $\pi_1(Y)=\pi_1 (\mathcal{K}_\varepsilon^+)=E/K$.

Hence the affirmative answer to the the D(2) problem would be equivalent to:

{\bf For all finitely presented groups $E$ and all perfect finitely closed normal subgroups $K\vartriangleleft E$ and all finite presentations $\varepsilon$ of $E$, there exists a finite presentation $\mathcal{G}$ of $E/K$ and a homotopy equivalence $\mathcal{K}_\varepsilon^+ \sim \mathcal{K}_\mathcal{G}$ inducing the identity $1\colon E/K \to E/K$ on fundamental groups.}

\lem{The following are equivalent: \newline \newline i) There exists a homotopy equivalence: $\mathcal{K}_\varepsilon^+ \sim \mathcal{K}_\mathcal{G}$ inducing the identity $1\colon E/K \to E/K$ on fundamental groups.\newline
\newline ii) There exists a chain homotopy equivalence: $C_*(\mathcal{K}_\varepsilon^+) \sim C_*(\mathcal{K}_\mathcal{G})$ over $\Z[E/K]$.}

\proof{i) $\Rightarrow$ ii) is immediate.  Conversely, from ii) we have a chain homotopy equivalence between the algebraic complexes associated to a finite cohomologically 2--dimensional 3--complex and a finite 2--complex (with respect to an isomorphism of fundamental groups).  To show that ii)$\Rightarrow$ i) we must construct a homotopy equivalence between the spaces, inducing the same isomorphism on fundamental groups.  For finite fundamental groups this is done in \cite[proof of Theorem 59.4]{John1}.  The same argument holds for all finitely presented fundamental groups \cite[Appendix B, {\sl Proof of Weak Realization Theorem}]{John1}.}

From the proof of lemma 2.4, $C_*(\mathcal{K}_\varepsilon)\otimes_E \Z[E/K] \sim C_*(\mathcal{K}_\varepsilon^+)$.  Hence we have:

\thm{The following two statements are equivalent: \newline \vspace{-3mm}\newline
i) Let $X$, a finite 3--complex with $H^3(X;\beta)=0$ for all coefficient bundles $\beta$.  Then $X$ is homotopy equivalent to a finite 2--complex.
\newline \vspace{-3mm}\newline
ii) Let $K$ be a perfect finitely closed normal subgroup of a finitely presented group $E$.  For each finite presentation $\varepsilon$ of $E$, there exists a finite presentation $\mathcal{G}$ of $E/K$, such that we have a chain homotopy equivalence over $\Z[E/K]$: $$C_*(\mathcal{K}_\varepsilon)\otimes_E \Z[E/K] \to C_*(\mathcal{K}_\mathcal{G})$$}

Suppose we have a short exact sequence: $$1\to L\to F\to G\to 1$$ where $G$ is a finitely presented group and $F$ is a free group generated by elements $g_1,\cdots,g_n$.  Let $R_1,\cdots,R_m$ be elements of $L$.

\define{$\langle g_1,\cdots,g_n \vert R_1, \cdots, R_m \rangle$ is called a finite partial presentation for $G$ when the normal closure $N_F(R_1,\cdots,R_m)$ surjects onto $L/[L,L]$ under the quotient map $L \to L/[L,L]$.}

Note that a finite partial presentation $\varepsilon=\langle g_1,\cdots,g_n \vert R_1, \cdots R_m \rangle$ as above is an actual finite  presentation of some group $E$, so it has a well defined Cayley complex $\mathcal{K}_\varepsilon$.  

Let $K$ denote the kernel of the homomorphism $E \to G$ sending each $g_i$ to the corresponding element in $G$.  If $G$ is finitely presented then it is finitely presented on the generators in $\varepsilon$ \cite[Chapter 1, Proposition 17]{Cohe}.  As $K$ is the normal closure in $E$ of the images of this finite set of relators we have that $K$ is finitely closed.  

Further $K$ is perfect as every $k \in K$ may be lifted to an element of $L$ which may be written in the form $ab$ where $a \in [L,L]$ and $b \in N_F(R_1,\cdots R_m)$.  Thus $k$ is equal to the image of $a$ in $E$, so $k \in [K,K]$.  Thus a finite partial presentation $\varepsilon$ of a finitely presented group $G$ may be viewed as a presentation satisfying the hypothesis' of statement ii) in Theorem 4.2.

Conversely, given $\varepsilon$ as in statement ii) Theorem 4.2, we have that $\varepsilon$ is a finite partial presentation of $E/K$ (as $K= [K,K]$), and $E/K$ is finitely presented (as K is finitely closed).  

Thus statement ii) is equivalent to:

ii)'  Given a finite partial presentation $\varepsilon$ of a finitely presented group $G$, there exists a finite presentation $\mathcal{G}$ of $G$, such that we have a chain homotopy equivalence:
$$C_*(\mathcal{K}_\varepsilon)\otimes_E \Z[G] \to C_*(\mathcal{K}_\mathcal{G})$$

where $E$ is the group presented by $\varepsilon$ and each $x \in E$ acts on $\ZG$ by left multiplication by its image in $G$.

One approach to the D(2) problem is to use Euler characteristic as an obstruction.  That is, given a finite cohomologically 2--dimensional 3--complex $X$, if we can show that every finite 2--complex $Y$ with $\pi_1(Y)=\pi_1(X)$ satisfies $\chi(X) < \chi(Y)$ then clearly $X$ cannot be homotopy equivalent to any such $Y$.  It has been shown that certain constructions involving presentations of a group would allow one to construct such a space \cite[Theorem 3.5]{Harl}.  A candidate for such a space is given in \cite{Brid}.  In light of Corollary 2.5 and Theorem 3.4 we are able to generalize this approach.

The deficiency ${\rm Def}(\mathcal{G})$ of a finite presentation $\mathcal{G}$ is the number of generators minus the number of relators.  We say a presentation of a group is minimal if it has the maximal possible deficiency.  A finitely presented group $G$ always has a minimal presentation, because an upper bound for the deficiency of a presentation is given by $\Rk(G /[G,G])$.  The deficiency ${\rm Def}(G)$ of a finitely presented group $G$ is defined to be the deficiency of a minimal presentation.

 Again let $K \vartriangleleft E$ be a perfect finitely closed normal subgroup.  Then if $\varepsilon$ is a finite presentation of $E$ and $\mathcal{G}$ is a finite presentation for $E/K$ we have: $$\chi(\mathcal{K}_\varepsilon^+) = \chi(\mathcal{K}_\varepsilon) = 1- {\rm Def}(\varepsilon)\qquad \qquad\qquad \chi(\mathcal{K}_\mathcal{G}) = 1-{\rm Def}(\mathcal{G})$$

 \lem{If ${\rm Def}(E)>{\rm Def}(E/K)$ then given a minimal presentation $\varepsilon$ of $E$ we have that $\chi(\mathcal{K}_\varepsilon^+)< \chi(\mathcal{K}_\mathcal{G})$ for any finite presentation $\mathcal{G}$ of $E/K$.}

  \hspace{2mm} Proof:  
  
  $\chi( \mathcal{K}_\mathcal{G})= 1-{\rm Def}(\mathcal{G}) \geq 1- {\rm Def}(E/K) > 1- {\rm Def}(E) = 1- {\rm Def}(\varepsilon) =\chi(\mathcal{K}_\varepsilon^+ )$ \hfill $\Box$

\bigskip
Suppose we have a short exact sequence of groups: $$1\to K\to E \to G \to 1$$ with $E$, $G$ finitely presented.  Then given a finite presentation for $E$, the images in $G$ of the generators will generate $G$.  We may present $G$ on these generators with a finite set of relators \cite[Chapter 1, Proposition 17]{Cohe}.  Let $k_1,\cdots,k_r$ denote the elements of $K$ represented by these relators.  Then $K$ is the normal closure in $E$ of $k_1,\cdots,k_r$ and so $K$ is finitely closed in $E$.  In particular $K/[K,K]$ is generated by the $k_1,\cdots,k_r$ as a right module over $\ZG$ (where $G$ acts on $K/[K,K]$ by conjugation).  Let ${\rm rk}_G(K)$ denote the minimal number of elements required to generate $K/[K,K]$ over $\ZG$.

\thm{The following statements are equivalent:\newline \newline
i) There exists a connected finite cohomologically 2--dimensional 3--complex $X$, such that for all finite connected 2--complexes $Y$ with $\pi_1(Y)=\pi_1(X)$ we have $\chi(X)<\chi(Y)$.\newline\newline ii) There exists a short exact sequence of groups $1\to K\to E \to G \to 1$ with $E$, $G$ finitely presented and $$ {\rm rk}_G(K) + {\rm Def}(G) < {\rm Def}(E)$$}

\hspace{2mm} Proof:  { i)$\Rightarrow$ii): By Theorem 3.4, $X$ is homotopy equivalent to $\mathcal{K}_\varepsilon^+$ for some finite presentation $\varepsilon$ of some group $E$ and some perfect finitely closed normal subgroup $K$.  Let $G= E/K$.  We have a short exact sequence: $$1 \to K \to E \to G \to 1$$  As $K$ is finitely closed, $G$ is finitely presented.  As $K$ is perfect we have ${\rm rk}_G(K) =0$.  Let $\mathcal{G}$ be some finite presentation of $G$.  We have: $$1-{\rm Def}(\varepsilon) = \chi(\mathcal{K}_\varepsilon^+) <\chi(\mathcal{K}_\mathcal{G}) = 1- {\rm Def}(\mathcal{G})$$  Thus ${\rm Def}(\mathcal{G}) < {\rm Def}(\varepsilon)$.   As $\mathcal{G}$ was chosen arbitrarily, we have ${\rm Def}(G) < {\rm Def}(\varepsilon) \leq {\rm Def}(E)$.  Hence ${0+ \rm Def}(G) < {\rm Def}(E)$ as required.
\newline\newline ${}$
\hspace{14mm}ii)$\Rightarrow$i):  We start with the short exact sequence $1\to K \to E \to G\to 1$.  Let $k_1,\cdots,k_r \in K$ generate $K/[K,K]$ over $\ZG$, where $r={\rm rk}_G(K)$.  Let $K'$ denote the normal closure in $E$ of $k_1,\cdots,k_r$.  Then we have a short exact sequence:$$1 \to K/K' \to E/K' \to G\to 1$$  Then $K=K'[K,K]$ so $K/K'$ is perfect.  From the discussion preceding this theorem we know that $K$ is finitely closed in $E$, so $K/K'$ must be finitely closed in $E/K'$.   Also $E/K'$ may be presented by taking a minimal presentation of $E$ and adding $r$ relators (representing to $k_1,\cdots,k_r$). Hence: $${\rm Def}(E/K') \geq {\rm Def}(E) - {\rm rk}_G(K) > {\rm Def}(G)$$ Take a minimal presentation $\varepsilon$ of $E/K'$ and let $X=\mathcal{K}_\varepsilon^+$, where $+$ is taken with respect to $K/K'$.  Any finite connected 2--complex $Y$ with $\pi_1(Y) = \pi_1(X)$ is homotopy equivalent to $\mathcal{K}_\mathcal{G}$ for some finite presentation $\mathcal{G}$ of $G$.  Therefore by lemma 4.4 we have $\chi(X) < \chi(Y)$ as required.
} \hfill $\Box$

\vspace{1.4mm}
We note that Michael Dyer proved ii)$\Rightarrow$i) in the case where $H^3(G;\ZG) =0$ and $E$ is a free group whose generators are the generating set for some minimal presentation of $G$ \cite[Theorem 3.5]{Harl}.

\bigskip
\noindent    Address: Building 54\\
    School of Mathematics\\
    University of Southampton\\
    University Road\\
    Southampton SO17 1BJ\\
    Telephone: +44 (0)2380 592 141\\
email:\verb|wajid@mannan.info|

\end{document}